\documentclass[11pt]{article}
\usepackage[]{array,theorem,amssymb,url}
\newtheorem{Theorem}{Theorem}
\newtheorem{Proposition}{Proposition}
\newtheorem{Corollary}{Corollary}
\newtheorem{Definition}{Definition}
\newtheorem{Example}{Example}
\newtheorem{Lemma}{Lemma}
\newtheorem{Remark}{Remark}
\newcommand{\lm}{{\mathbb{F}}_q}

\newcommand{\firkant}{\begin{flushright} \vspace{-.7cm} $\Box$ \end{flushright}}\newcommand{\funcf}{{\mathbb{F}}/{\mathbb{F}}_q}

\newcommand{\Lamdba}{\Lambda}
\newcommand{\lmto}{{\mathbb{F}}_{q^2}}

\newcommand{\PP}{{\mathcal{P}}}
\begin{document}
\begin{center}
{\Large{Bounding the number of rational places\\
using
    Weierstrass semigroups}}\\
\ \\
Olav Geil\\
Department of Mathematical Sciences\\
Aalborg University\\
9220 Denmark\\

\ \\

Ryutaroh Matsumoto\\
Department of Communications and Integrated Systems\\
Tokyo Institute of Technology\\
152-8550 Japan\\

\ \\

Email: olav@math.aau.dk and ryutaroh@rmatsumoto.org
\end{center}
\ \\
\begin{center}
\begin{minipage}{12cm}
{\small{{\textbf{Abstract:}} Let $\Lambda$ be a numerical
    semigroup. Assume there exists an algebraic function field over
    ${\mathbb{F}}_q$ in one variable which possesses a rational place
    that has $\Lambda$ as its Weierstrass semigroup. We ask the question
    as to how many rational places such a function field can possibly have
    and we derive an upper bound in terms of the generators of $\Lambda$
    and $q$. Our bound is an improvement to Lewittes' bound
    in~\cite{lewittes} which takes into account only the multiplicity
    of $\Lambda$ and $q$. From the new bound we derive significant
    improvements to Serre's upper bound in the cases $q=2$, $3$ and
    $4$. We finally show that Lewittes' bound has important
    implications to the theory of towers of function fields.\\  


\noindent {\textbf{MSC-class:}}  11G20; 14G15 (primary); 14H05; 14H25}}
\end{minipage}
\end{center}
\ \\

\section{Introduction}\label{secintro}
Throughout this paper by a function field we will always
mean an algebraic function field of one variable unless otherwise
stated. Given a function field $\funcf$ we denote by
$N({\mathbb{F}})$ the number of rational places and we denote by
$g({\mathbb{F}})$ the genus. We will always assume that
${\mathbb{F}}_q$ is the full constant field of ${\mathbb{F}}$. For applications in coding theory it is
desirable to have $N({\mathbb{F}})/g({\mathbb{F}})$ as high as
possible as this allows for the construction of codes with good
parameters. The above observation has led to an extensive research on
the problem of deciding given a fixed constant field ${\mathbb{F}}_q$
and a fixed number $g$ what is the highest number $N_q(g)$ such that
a 
function field $\funcf$ exists with $N({\mathbb{F}})=N_q(g)$ and
$g({\mathbb{F}})=g$.\\
Recall that, for any rational place the number of gaps in the corresponding
Weierstrass semigroup $\Lambda$ equals the genus $g$ of the
corresponding function field. That is, $\# ({\mathbb{N}}_0
\backslash \Lambda )=g$.  This suggests that in some cases a Weierstrass semigroup
$\Lambda$ for a rational place might hold more information about the
number of rational places of the function field than does the genus
alone. This theme was firstly explored by Lewittes
in~\cite{lewittes}, though the bound by St\"{o}hr and Voloch (\cite[pp.\ 14-15]{storvoloch}) induces a
bound in terms of a Weierstrass semigroup under certain conditions. The smallest non zero
element in a numerical semigroup $\Lambda$ is called the multiplicity of $\Lambda$ and we
denote it by $\lambda_1$. Using function field theory
Lewittes showed that if $\lambda_1$ is the multiplicity of a
Weierstrass semigroup corresponding to a rational place of
${\mathbb{F}}/{\mathbb{F}}_q$ then $N({\mathbb{F}}) \leq q
\lambda_1+1$ holds. In the present paper we derive an improved upper bound on
$N({\mathbb{F}})$ as we take into account not only the multiplicity
but also all  the other elements in a generating set of
$\Lambda$. Let 
$$\Lambda=\langle \lambda_1, \ldots , \lambda_m
\rangle=\{a_1\lambda_1+\cdots a_m\lambda_m \mid a_1, \ldots ,a_m \in {\mathbb{N}}_0\}$$
$0 < \lambda_1 < \cdots < \lambda_m$ be the Weierstrass semigroup of a rational place in a
function field $\funcf$. Our new finding is that $N({\mathbb{F}})\leq \#
(\Lambda \backslash \cup_{i=1}^{m} (q\lambda_i + \Lambda ))+1$
holds. Here, $\gamma +\Lambda$ means $\{\gamma+\lambda \mid \lambda \in \Lambda
\}$. Lewittes' bound can be viewed as a corollary to the new bound as  $\# (\Lambda \backslash (q
\lambda_1+\Lambda))+1=q \lambda_1 +1$ holds.
The new bound is often 
better than Lewittes' bound. For $q$ being $2$, $3$ and $4$ we get as a
corollary to the new bound a significant improvement to the Serre
bound. As will be demonstrated in this paper Lewittes' bound has
implications to the theory of towers of function fields. We show that Lewittes' bound holds
information that describes rather accurately certain aspects of the second
tower of Garcia and Stichtenoth. Finally,  we show that one cannot hope to
construct asymptotically good towers of function fields having
telescopic Weierstrass semigroups.\\
The paper is organized as follows. In Section~\ref{secproof} we prove
the bounds and state some corollaries.  The rest of the paper is devoted to
investigating how good are the bounds and how
much does the new bound improve upon Lewittes' bound. In
Section~\ref{secg8} we first deal with a selection of semigroups of
high genus and we then apply for $q$ being equal to $2$, $3$ and $4$ the bounds to all semigroups of genus
$8$. In
Section~\ref{secnq} we estimate $N_q(g)$ for $q$ being equal to $2$, $3$ and $4$.
Finally, in
Section~\ref{secasymp} we are concerned with asymptotically good towers
of function fields. Section~\ref{secon} is the conclusion.
\section{The bounds}\label{secproof}
Throughout this paper let $\Lambda$ be a numerical semigroup with
finitely many gaps (meaning that ${\mathbb{N}}_0 \backslash \Lambda$
is finite) and let $\{ \lambda_1, \ldots ,
\lambda_m\}$ be a generating set for $\Lambda$ with  $0 <\lambda_1 <
\cdots < \lambda_m$. The reader may think of
$\lambda_1, \ldots , \lambda_m$ as being a minimal generating set but
being so or not is actually of no implication. 
\begin{Definition}
Let $\Lambda$ be a numerical semigroup with finitely many gaps. If there
does not exist a function field over ${\mathbb{F}}_q$ having a
rational place which Weierstrass semigroup equals $\Lambda$ then we write
$N_q(\Lambda)=0$. If such function fields exist we define
\begin{eqnarray}
N_q(\Lambda)&=&\max\{N({\mathbb{F}}) \mid {\mathbb{F}} {\mbox{ \  is a
    function field over \ }} {\mathbb{F}}_q {\mbox{\ having a
    rational}} \nonumber \\
&&{\mbox{ \ \ \ \ \ \ \ \ \ \ \ \ \ \ \ \ \ \ \ \ \ \ \ \ place which Weierstrass semigroup equals \ }} \Lambda \}.\nonumber
\end{eqnarray}
\end{Definition}
Our main result is:
\begin{Theorem}\label{themain}
Let $\Lambda=\langle \lambda_1, \ldots , \lambda_m \rangle$ be a
numerical semigroup with finitely many gaps, $0 < \lambda_1 < \cdots
<\lambda_m$. We have
\begin{eqnarray}
N_q(\Lambda)&\leq &\# (\Lambda \backslash \cup_{i=1}^{m} (q\lambda_i + \Lambda ))+1 \label{stjerneskud2}
\end{eqnarray}
which implies 
\begin{eqnarray}
N_q(\Lambda)& \leq & \# (\Lambda \backslash (q
\lambda_1+\Lambda))=q \lambda_1 +1. \label{stjerneskud3}
\end{eqnarray}
Here, $\gamma +\Lambda$ means $\{\gamma+\lambda \mid \lambda \in \Lambda
\}$. 
\end{Theorem}
\noindent {\it{Proof:}} \ \   Let $\funcf$ be a function field. Let its rational places be ${\mathcal{P}}_1, \ldots ,
{\mathcal{P}}_{N-1}, {\mathcal{P}}$ and assume that the Weierstrass
semigroup corresponding to ${\mathcal{P}}$ is $\Lambda$. 
Define $R=\cup_{s=0}^{\infty}
{\mathcal{L}}(s{\mathcal{P}})$ and let $R_{-1}=\{0\}$ and
$R_t=\cup_{s=0}^t {\mathcal{L}}(s{\mathcal{P}})$ for $t \in
{\mathbb{N}}_0$. It is well known that
\begin{eqnarray}
\begin{array}{ll}
R_t=R_{t-1} & {\mbox{ \ if \ }} t \in {\mathbb{N}}_0 \backslash \Lambda\\
\dim(R_t)=\dim(R_{t-1})+1 & {\mbox { \ if }} t \in \Lambda.
\end{array} \label{eqsnabel1}
\end{eqnarray}
Here $\dim$ denotes the dimension as a vector space over
${\mathbb{F}}_q$. Let $\varphi : R \rightarrow {\mathbb{F}}_q^{N-1}$ be the map
$\varphi(f)=(f({\mathcal{P}}_1), \ldots , f({\mathcal{P}}_{N-1}))$ and define
$E_t=\varphi(R_t)$ for $t \in {\mathbb{N}}_0 \cup
\{-1\}$. From~(\ref{eqsnabel1}) we observe that $\dim(E_{-1})=0$ and that
$\dim(E_t)=\dim(E_{t-1})$ for all $t \in {\mathbb{N}}_0 \backslash \Lambda$. For $t \in
\Lambda$ we can either have $\dim(E_t)=\dim(E_{t-1})$ or 
$\dim(E_t)=\dim(E_{t-1})+1$. The map $\varphi$ is surjective meaning
that for $t$ large enough $\dim(E_t)=N-1$. Hence, if we can give an
  upper bound on the number of $t \in \Lambda$ for which
  $\dim(E_t)=\dim(E_t-1)+1$ holds then this upper bound will also be
  an upper bound on the number $N-1$. To prove~(\ref{stjerneskud2})
  we therefore only need to show that $\dim(E_t)=\dim(E_{t-1})+1$
  cannot happen when $t \in q\lambda_i+\Lambda$ for some $i$. For this
  purpose let for $i=1, \ldots , m$,  $x_i \in R$ be an element with
  $-v_{{\mathcal{P}}}(x_i)=\lambda_i$. Here $v_{{\mathcal{P}}}$ is the
  valuation corresponding to ${\mathcal{P}}$. Given $t=q \lambda_i+\lambda$
  with $\lambda \in \Lambda$ choose $f \in R_{\lambda} \backslash
  R_{\lambda-1}$. We have $x_i^qf \in R_t \backslash R_{t-1}$ and
  $x_if \in R_{t-1}$. Clearly, $\varphi(x_i^qf)=\varphi(x_if)$ and
  the
  proof of~(\ref{stjerneskud2}) is complete. Finally to see
  that~(\ref{stjerneskud3}) follows from~(\ref{stjerneskud2}) we
  note that $\# (\Lambda \backslash (\cup_{i=1}^m q \lambda_i
  +\Lambda)) \leq \# (\Lambda \backslash q \lambda_1
  +\Lambda)$ and we use the following lemma. \\
\ \ \ \ \firkant 
\begin{Lemma}\label{lemfornum}
Let $\Lambda \subseteq {\mathbb{N}}_0$ be a semigroup with finitely
many gaps. For
  $\lambda \in \Lambda$ we have $\# (\Lamdba \backslash (\lambda +
  \Lambda))= \lambda$.
\end{Lemma}
\noindent {\it{Proof:}} \ \ See~\cite[Lem.\ 5.15]{handbook}. \firkant 
\begin{Example}
In this example we apply Theorem~\ref{themain} to the semigroup
$$\Lambda=\langle \lambda_1=3,\lambda_2=5\rangle$$ 
in the
case $q=2$. We have
$$\Lambda \backslash (q\lambda_1+\Lambda)= \{ 0,3, 5, 8, 10, 13\},$$
whereas
$$\Lambda \backslash ((q\lambda_1+\Lambda)\cup (q\lambda_2+\Lambda))=
\{0,3, 5, 8\}.$$ Lewittes' bound~(\ref{stjerneskud3}) states
$N_2(\Lambda) \leq 7$, whereas the new bound~(\ref{stjerneskud2}) gives $N_2(\Lambda) \leq 5$.
\end{Example}
The following Proposition gives us some information on how good or bad
the bound in~(\ref{stjerneskud2}) can possibly be.
\begin{Proposition}\label{proboundbound}
We have $$q \lambda_1+1-g \leq \# (\Lambda \backslash ( \cup_{i=1}^m
(q\lambda_i+\Lambda)))+1 \leq \min\{q \lambda_1+1,  q^m+1\}.$$
\end{Proposition}
\noindent {\it{Proof:}} \ \ To see the first inequality observe that there are at
least $q \lambda_1-g$ elements in $\Lambda$ that are smaller than
$q\lambda_1$. Regarding the last inequality the upper bound $\lambda_1
q+1$ comes from Theorem~\ref{themain}. To see the upper bound $q^m+1$ we
note that all $\lambda \in \Lambda$ can be written as
$a_1\lambda_1+ \cdots + a_m\lambda_m$ for some $a_1, \ldots , a_m \in
{\mathbb{N}}_0$. If $\lambda \in \Lambda \backslash
(\cup_{i=1}^m(q\lambda_i+\Lambda)$ then necessarily $a_1, \ldots , a_m
  < q$ must hold. \firkant
Recall, that by $N_q(g)$ we denote the maximal number of rational
places that a function field can possibly have if its full constant
field is ${\mathbb{F}}_q$ and its genus is $g$. Serre's bound
\begin{eqnarray}
\mid N_q(g)-(g+1)\mid & \leq & g \lfloor 2 \sqrt{q} \rfloor \nonumber
\end{eqnarray}
implies that if $\Lambda$ is of genus $g$ then
\begin{eqnarray}
N_q(\Lambda) &\leq&g \lfloor 2 \sqrt{q} \rfloor +q + 1 \label{eqsnabel12}
\end{eqnarray}
holds. Writing $r=\frac{\lambda_1-1}{g}$ we see that for
$$r\leq \frac{\lfloor 2\sqrt{q}\rfloor}{q}$$
the right side
of~(\ref{stjerneskud3}) (and therefore also the right side
of~(\ref{stjerneskud2})) is always better than~(\ref{eqsnabel12}).
On the contrary Proposition~\ref{proboundbound} tells us that the new
bound~(\ref{stjerneskud2}) can not produce a number smaller than
$q\lambda_1+1-g$. We conclude that for 
\begin{eqnarray}
r& >&\frac{\lfloor 2\sqrt{q} \rfloor +1}{q} \label{eqsnabel13}
\end{eqnarray}
the bound~(\ref{stjerneskud2}) can not compete with~(\ref{eqsnabel13}). Observe however, that as $r$ is at
most $1$ the condition in~(\ref{eqsnabel13}) is never
satisfied for $q=2,3,4,5$.\\
We now present some corollaries to
Theorem~\ref{themain}. The first corollary is a trivial restatement of
Lewittes' bound~(\ref{stjerneskud3}).
\begin{Corollary}\label{corrat1}
If $\funcf$ possesses $N$ rational places then for all
of the corresponding 
Weierstrass semigroups we have 
\begin{eqnarray}
\lambda_1 \geq (N-1)/q. \label{eqsemsoc}
\end{eqnarray}
\end{Corollary}
\begin{Example}
The norm-trace function field defined by $x^{{q^r-1}/(q-1)} -
y^{q^{r-1}}-y^{q^{r-2}}- \cdots - y$ has $N=q^{2r-1}+1$ rational
places (see~\cite{normtrace}). All but one correspond to affine points of the curve. The
last one is denoted by $P_{\infty}$. It is well-known that the
Weierstrass semigroup of $P_{\infty}$ equals $\langle q^{r-1},
(q^r-1)/(q-1) \rangle $. That is, $\lambda_1 = q^{r-1}$.  But also
$(N-1)/q^r$ equals $q^{r-1}$ and the norm-trace function field
therefore is an example where the bound in Corollary~\ref{corrat1} is reached.
\end{Example}
\begin{Corollary}\label{cort}
Define 
$$t=\# \{ \lambda\in \Lambda \mid \lambda \in [ \lambda_1+1,
\lambda_1+\lceil \lambda_1/q \rceil -1 ] \}.$$
We have
\begin{equation}
N_q(\Lambda) \leq q \lambda_1 -t+1. \label{eqo1}
\end{equation}
\end{Corollary}
\noindent {\it {Proof:}} \ \ For $\lambda \in \Lambda$ with $\lambda \in  [ \lambda_1+1,
\lambda_1+\lceil \lambda_1/q\rceil -1 ] $ we have $q \lambda \neq q \lambda_1 +
\eta$ for any $\eta \in \Lambda$ as there are no non zero $\eta \in
\Lambda$ with $\eta < \lambda_1$. Therefore the number
on the right side of~(\ref{stjerneskud2}) is at least $t$ smaller than
the number on the right side of~(\ref{stjerneskud3}). \firkant
\begin{Example}\label{exordinary}
Consider the
case $\lambda_1=g+1$. That is, the case $\Lambda=\{0,g+1, g+2, \ldots
\}$. The number $t$ from Corollary~\ref{cort} becomes equal to $\lceil
(g+1)/q\rceil -1$. Hence 
\begin{equation}
N_q(\Lambda ) \leq q(g+1)+2-\lceil (g+1)/q \rceil \label{eqordinary}
\end{equation}
holds. Given $\lambda > \lambda_1$ we have $q\lambda \notin q\lambda_1+\Lambda$ if and only if $\lambda \in [
\lambda_1+1,\lambda_1+\lceil \lambda_1 / q \rceil -1]$ and $q\lambda +
\eta \in q\lambda_1+\Lambda$ holds for all $\eta \in \Lambda
\backslash \{ 0 \}$. Hence, for the particular semigroup in the
present example, we have 
$$\#(\Lambda
\backslash(\cup_{i=1}^m(\lambda_i+\Lambda)))= q \lambda_1
-t+1=q(g+1)+2-\lceil (g+1)/q \rceil .$$
\end{Example}
The following remark gives some criteria under which the
bounds~(\ref{stjerneskud2}) and (\ref{stjerneskud3}) are the same.
\begin{Remark}\label{remthesame}\label{remone}
The conductor of a semigroup $\Lambda \subseteq {\mathbb{N}}_0$ with
finitely many gaps is the smallest number $c$ such that there are no
gaps greater or equal to $c$. 
The conductor is known to be smaller or equal
to 
$2g$ (\cite[Pro.\ 5.7]{handbook}). If $q \lambda_1 + c \leq q \lambda_2$ then it is
clear that the number on the right side of~(\ref{stjerneskud2}) is the
same as the number on the right side 
of~(\ref{stjerneskud3}). In particular the numbers are the
same if $q \lambda_1 + 2g \leq q \lambda_2$ 
\end{Remark}
We conclude this section by mentioning that by
using the theory of algebraic geometric codes one  can
sharpen~(\ref{stjerneskud2}) to 
\begin{equation}
N_q(\Lambda)\leq \# \left(\left( \Lambda \cap \{\lambda \mid \lambda \leq
  N'+2g-2\}\right) \backslash (\cup_{i=1}^{m}(q\lambda_i+\Lambda))\right)+1. \nonumber
\end{equation}
Here, $N'$ is some a priori known upper bound on $N_q(\Lambda)$, e.\
g.\ $N'=N_q(g)$. We should
mention that in our experiments we have not been able to find any
example where this improvement gives a number that is simultaneously
smaller than both $N'$ and the right side of~(\ref{stjerneskud2}). \\
In the following sections we investigate how good is the new
bound~(\ref{stjerneskud2}) and in particular how good is it compared
to Lewittes' bound. We start by investigating a selection of concrete
semigroups. 
\section{Examples }\label{secg8}
In this section we apply the bounds~(\ref{stjerneskud2}) and
  (\ref{stjerneskud3}) to a number of concrete semigroups. 
\begin{Example}\label{exgstor}
In Table~\ref{tabmangeq} we consider a collection of $7$
semigroups. We apply the bounds to a number of fields of characteristics $2$ and
$3$. Restricting to characteristics $2$ and $3$ allows us to get
information on the number $N_q(g)$ from van der Geer and van
der Vlugt's table in~\cite{geervlugt}. An entry $x/y$ in the row named
``bounds'' indicates that Lewittes' bound
produces $x$ and that the new bound produces $y$. An interval in the
row named $N_q(g)$ means that $N_q(g)$ is known to be in this
interval. An $\ast$ in the same row means that the table
in~\cite{geervlugt} is empty. Table~\ref{tabmangeq} illustrates that
the new bound can be quite an improvement to Lewittes' bound and that
it can be much smaller than $N_q(g)$ also when Lewittes' bound is
not. It is clear that we get the most significant results for smallest $q$.
\begin{table}
{\small{
$$
\begin{array}{|c|}
\hline
\Lambda=\langle 8,9,20 \rangle {\mbox{ \ \ \ \ \ \ \ }} g=20\\
\hline
\begin{array}{c||c|c|c|c|c|c} 
q&2&3&4&8&9&16\\
\hline
{\mbox{bounds}}&17/9&25/16&33/25&65/65&73/73&129/129\\
\hline
N_q(g)&19-21&30-34&40-45&76-83&70-91&127-139\\
\end{array}\\
\hline
\end{array}
$$}}
{\small{
$$
\begin{array}{|c|}
\hline
\Lambda=\langle 13,15,17,18,20 \rangle {\mbox{ \ \ \ \ \ \ \ }} g=23\\
\hline
\begin{array}{c||c|c|c|c|c|c} 
q&2&3&4&8&9&16\\
\hline
{\mbox{bounds}}&27/14&40/30&53/46&105/102&118/118&209/209\\
\hline
N_q(g)&22-23&32-37&45-50&68-92&92-101&126-155\\
\end{array}\\
\hline
\end{array}
$$}}
{\small{
$$
\begin{array}{|c|}
\hline
\Lambda=\langle 13,15,24,31 \rangle {\mbox{ \ \ \ \ \ \ \ }} g=38\\
\hline
\begin{array}{c||c|c|c|c|c|c} 
q&2&3&4&8&9&16\\
\hline
{\mbox{bounds}}&27/13&40/28&53/40&105/97&118/112&209/207\\
\hline
N_q(g)&30-33&\ast&64-74&129-135&105-149&193-233\\
\end{array}\\
\hline
\end{array}
$$}}
{\small{
$$
\begin{array}{|c|}
\hline
\Lambda=\langle 20,22,23,24,26 \rangle {\mbox{ \ \ \ \ \ \ \ }} g=47\\
\hline
\begin{array}{c||c|c|c|c|c|c} 
q&2&3&4&8&9&16\\
\hline
{\mbox{bounds}}&41/15&61/34&81/57&161/147&181/166&321/313\\
\hline
N_q(g)&36-38&54-65&73-87&126-161&154-177&\ast\\
\end{array}\\
\hline
\end{array}
$$}}
{\small{
$$
\begin{array}{|c|}
\hline
\Lambda=\langle 13,14,20 \rangle {\mbox{ \ \ \ \ \ \ \ }} g=42\\
\hline
\begin{array}{c||c|c|c|c|c|c} 
q&2&3&4&8&9&16\\
\hline
{\mbox{bounds}}&27/9&40/17&53/33&105/95&118/102&209/195\\
\hline
N_q(g)&33-35&52-59&75-80&129-147&122-161&209-254\\
\end{array}\\
\hline
\end{array}
$$}}
{\small{
$$
\begin{array}{|c|}
\hline
\Lambda=\langle 16,17,18,19 \rangle {\mbox{ \ \ \ \ \ \ \ }} g=45\\
\hline
\begin{array}{c||c|c|c|c|c|c} 
q&2&3&4&8&9&16\\
\hline
{\mbox{bounds}}&33/9&49/19&65/32&129/108&145/124&257/257\\
\hline
N_q(g)&33-37&54-62&80-84&144-156&136-170&242-268\\
\end{array}\\
\hline
\end{array}
$$}}
{\small{
$$
\begin{array}{|c|}
\hline
\Lambda=\langle 10,11,20,22 \rangle {\mbox{ \ \ \ \ \ \ \ }} g=45\\
\hline
\begin{array}{c||c|c|c|c|c|c} 
q&2&3&4&8&9&16\\
\hline
{\mbox{bounds}}&21/5&31/10&41/17&81/65&91/82&161/141\\
\hline
N_q(g)&33-37&54-62&80-84&144-156&136-170&242-268\\
\end{array}\\
\hline
\end{array}
$$}}
\caption{Semigroups from Example~\ref{exgstor}.}
\label{tabmangeq}
\end{table}
\end{Example}
\begin{Example}\label{eksg8}
From~\cite{rivaldo} we get all semigroups of genus $8$. There are $66$
of them. In Table~\ref{tabg8a} and Table~\ref{tabg8b} we then apply
the bounds~(\ref{stjerneskud2}) and (\ref{stjerneskud3}) to the cases
of $q$ being $2$, $3$ and $4$. As in the previous example an entry
$x/y$ means that Lewittes' bound produces $x$ whereas the new bound
produces $y$. From~\cite{geervlugt} we know that $N_2(8)=11$, $N_3(8)
\in \{17,18\}$ and $N_4(8) \in \{ 21, 22,23,24\}$. Lewittes' bound
tells us that in a function field over ${\mathbb{F}}_2$ of genus $8$
and with $N_2(8)=11$ rational places $13$ semigroups are not allowed
as Weierstrass semigroups of a rational place. The new bound gives us
that $33$ semigroups are not allowed. Assuming that $N_3(8)=18$
Lewittes' bound excludes $26$ semigroups whereas the new bound
excludes $31$ semigroups. Assuming $N_4(8)=24$ we get the exact same picture.
\begin{table}
{\small{
$$
\begin{array}{|l|c|c|c|}
\hline 
{\mbox{Semigroup}}&q=2&q=3 &q=4 \\
\hline 
\langle 2,17 \rangle &5/5&7/7&9/9\\
\langle 3,10,17 \rangle &7/6&10/10&13/13\\
\langle 3,11,16 \rangle &7/7&10/10&13/13\\
\langle 3,13,14 \rangle &7/7&10/10&13/13\\
\langle 4,6,13 \rangle &9/9&13/13&17/17\\
\langle 4,6,15,17 \rangle &9/9&13/13&17/17\\
\langle 4,7,17 \rangle &9/6&13/11&17/17\\
\langle 4,9,10 \rangle & 9/9&13/12&17/17\\
\langle 4,9,11 \rangle &9/7&13/13&17/17\\
\langle 4,9,14,15 \rangle &9/8&9/9&17/17\\
\langle 4,10,11,17 \rangle &9/9&13/13&17/17\\
\langle 4,10,13,15 \rangle &9/9&13/13&17/17\\
\langle 4, 11,13,14 \rangle &9/9&13/13&17/17\\
\langle 5,6,13 \rangle &11/7&16/12&21/18\\
\langle 5,6,14 \rangle &11/7&16/12&21/19\\
\langle 5,7,9 \rangle  &11/7&16/13&21/19\\
\langle 5,7,11 \rangle &11/9&16/14&21/19\\
\langle 5,7,13,16\rangle &11/8&16/14&21/20\\
\langle 5,8,9 \rangle  &11/9&16/15&21/20\\
\langle 5,8,11,12 \rangle &11/9&16/14&21/21\\
\langle 5,8,11,14,17 \rangle &11/9&16/15&21/20\\
\langle 5,8,12,14 \rangle  &11/9&16/15&21/21\\
\langle 5,9,11,12 \rangle &11/9&16/16&21/21\\
\langle 5,9,11,13,17 \rangle &11/9&16/15&21/21\\
\langle 5,9,12,13,16 \rangle  &11/10&16/16&21/21\\
\langle 5,11,12,13,14 \rangle &11/11&16/16&21/21\\
\hline
\end{array}
$$}}
\caption{Semigroups from Example~\ref{eksg8}}
\label{tabg8a}
\end{table}

\begin{table}
\label{tabone}
{\small{
$$
\begin{array}{|l|c|c|c|}
\hline 
{\mbox{Semigroup}}&q=2&q=3 &q=4 \\
\hline 
\langle 6,7,8,17 \rangle &13/8&19/15&25/22\\
\langle 6,7,9,17 \rangle  &13/10&19/17&25/22\\
\langle 6,7,10,11 \rangle &13/11&19/16&25/21\\
\langle 6,7,10,15\rangle &13/10&19/17&25/23\\
\langle 6,7,11,15,16 \rangle  &13/9&19/16&25/23\\
\langle 6,8,11,13,15 \rangle &13/11&19/19&25/25\\
\langle 6,8,10,13,15,17 \rangle &13/12&19/19&25/25\\
\langle 6,8,10,11,15 \rangle  &13/12&19/19&25/25\\
\langle 6,8,10,11,13 \rangle &13/11&19/18&25/25\\
\langle 6,8,9,10 \rangle &13/11&19/19&25/25\\
\langle 6,8,9,11 \rangle  &13/10&19/19&25/25\\
\langle 6,8,9,13 \rangle &13/11&19/19&25/25\\
\langle 6,9,10,11,14 \rangle &13/12&19/19&25/25\\
\langle 6,9,10,11,13 \rangle  &13/11&19/19&25/25\\
\langle 6,9,10,13,14,17 \rangle &13/12&19/19&25/25\\
\langle 6,9,11,13,14,16\rangle &13/12&19/19&25/25\\
\langle 6,10,11,13,14,15 \rangle  &13/12&19/19&25/25\\
\langle 7,8,9,10,11 \rangle &15/10&22/18&29/26\\
\langle 7,8,9,10,12 \rangle &15/10&22/18&29/26\\
\langle 7,8,9,10,13 \rangle  &15/10&22/18&29/26\\
\langle 7,8,9,11,12 \rangle &15/11&22/18&29/27\\
\langle 7,8,9,11,13 \rangle &15/11&22/18&28/27\\
\langle 7,8,9,12,13 \rangle  &15/11&22/18&29/27\\
\langle 7,8,10,12,13 \rangle &15/12&22/19&29/27\\
\langle 7,8,10,11,12\rangle &15/11&22/19&29/29\\
\langle 7,8,10,11,13 \rangle  &15/11&22/19&29/27\\
\langle 7,8,11,12,13,17 \rangle &15/12&22/20&29/28\\
\langle 7,9,10,11,12,13\rangle &15/11&22/20&29/27\\
\langle 7,9,10,11,13,15 \rangle  &15/11&22/20&29/28\\
\langle 7,9,10,12,13,15 \rangle &15/12&22/21&29/28\\
\langle 7,9,11,12,13,15,17 \rangle &15/12&22/21&29/28\\
\langle 7,10,11,12,13,15,16 \rangle  &15/13&22/21&29/29\\
\langle 8,9,10,11,12,13,14 \rangle &17/13&25/22&33/31\\
\langle 8,9,10,11,12,13,15 \rangle &17/13&25/22&33/31\\
\langle 8,9,10,11,12,14,15 \rangle  &17/13&25/22&33/31\\
\langle 8,9,10,11,13,14,15 \rangle &17/13&25/22&33/31\\
\langle 8,9,10,12,13,14,15 \rangle  &17/14&25/22&33/32\\
\langle 8,9,11,12,13,14,15 \rangle &17/14&25/24&33/32\\
\langle 8,10,11,12,13,14,15,17 \rangle &17/15&25/23&33/33\\
\langle 9,10,11,12,13,14,15,16,17 \rangle  &19/15&28/26&37/35\\
\hline
\end{array}
$$}}
\caption{Semigroups from Example~\ref{eksg8}}
\label{tabg8b}
\end{table}
\end{Example}
\section{Bounds on $N_q(g)$}\label{secnq}
From Lewittes' bound~(\ref{stjerneskud3}) we immediately get $N_q(g) \leq
q(g+1)+1$ as the multiplicity of a semigroup with $g$ gaps can be at
most $g+1$. This fact is not stressed in~\cite{lewittes} as the paper
contains slightly better bounds on $N_q(g)$ namely $N_q(g) \leq qg+2$
(\cite[Th.\ 1, part (a)]{lewittes}) and $N_2(g) \leq 2g-2$ (\cite[Eq.\
(19)]{lewittes}). The latter bounds are slightly better than Serre's upper bound
in the case of $q$ being equal to $2$, $3$ and $4$. 
We now investigate the implication of the new
result~(\ref{stjerneskud2}) for establishing bounds on
$N_q(g)$. We get the following proposition.
\begin{Proposition}\label{pronq}
\begin{equation}
N_q(g) \leq (q-\frac{1}{q})g+q+2-\frac{1}{q}. \label{eqo2} 
\end{equation}
\end{Proposition}
\noindent {\it{Proof:}} \ \ 
The proof uses Corollary~\ref{cort}. A 
reasonable estimate of the number $t$ in Corollary~\ref{cort} can be given in terms of $\lambda_1$ and $g$
alone. We have
\begin{eqnarray}
t&\geq&\max \{ \lceil \lambda_1/q \rceil - 1 -(g-(\lambda_1-1)),0\} \label{eqdoimand} \\
&\geq&\max \{\frac{\lambda_1}{q}+\lambda_1-g-2,0\} \label{eqdiamond1}
\end{eqnarray}
as there are $g-(\lambda_1-1)$ gaps greater than $\lambda_1$. 
Observe that $\frac{\lambda_1}{q}+\lambda_1-g-2 \leq 0$ holds for
$\lambda_1 \leq \frac{q}{q+1}(g+2)$. Hence, $N_q(g)\leq \max\{
K_1,K_2\}$ where 
\begin{eqnarray}
K_1&:=&q(\frac{q}{q+1}(g+2))+1 \nonumber \\
K_2&:=&\max\{q \lambda_1-(\frac{\lambda_1}{q}+\lambda_1-g-2)+1 \mid
\frac{q}{q+1}(g+2) \leq \lambda_1 \leq g+1\}. \nonumber
\end{eqnarray}
The maximal value of $q
\lambda_1-(\frac{\lambda_1}{q}+\lambda_1-g-2)+1$ is attained for
$\lambda_1=g+1$ and the proposition follows. \firkant
Observe, that the bound~(\ref{eqo2}) was obtained by showing
that the semigroup considered in Example~\ref{exordinary} is the worst
case. That is, (\ref{eqo2}) is almost the same as (\ref{eqordinary}). With the last part of Example~\ref{exordinary} in mind we cannot hope to
improve upon Proposition~\ref{pronq} using our method.\\
For $q$ being equal to $2$, $3$
and $4$ Proposition~\ref{pronq} is much better than Serre's upper 
bound. We get 
\begin{eqnarray}
N_2(g) & \leq &1 \frac{1}{2}g+3\frac{1}{2}, \nonumber \\
N_3(g) & \leq &2 \frac{2}{3}g+4\frac{2}{3}, \nonumber \\
N_4(g) & \leq &3 \frac{3}{4}g+5\frac{3}{4}, \nonumber 
\end{eqnarray}
whereas Serre's upper bound states
\begin{eqnarray}
N_2(g) & \leq & 2g+3, \nonumber \\
N_3(g) & \leq &3g+4, \nonumber \\
N_4(g) & \leq &4g+5. \nonumber 
\end{eqnarray}
For values of $q$ greater than or equal to $5$ Serre's upper bound is 
much better than Proposition~\ref{pronq}.
Even though Proposition~\ref{pronq} is better than Serre's upper bound for
three values of $q$ it does not provide 
information that is not already known. For $q$ being equal to $2$, $3$ and $4$
and $g$ being not too small Ihara's bound
\begin{eqnarray}
N_q(g) & \leq &q+1+ \lfloor (\sqrt{(8g+1)g^2+4(q^2-q)g}-g)/2 \rfloor .
\end{eqnarray}
namely outperforms the bound in Proposition~\ref{pronq}. For $q=2$ of course
the bound 
$$N_2(g) \leq (0.83)g+5.35,$$
which has been produced by the Oesterlé-Serre method
(see~\cite[Ex.\ 1.6.19]{oesterleserre}) even more outperforms Proposition~\ref{pronq}.
\section{Towers of  function fields}\label{secasymp}
With the results in the previous sections in mind unsurprisingly 
Lewittes' bound has implications for the theory of asymptotically good
towers of function fields.
Recall, that a sequence of function fields $(F^{(1)}/\lm, F^{(2)}/\lm \cdots )$ is called a
tower if $F^{(i)} \subseteq F^{(i+1)}$ holds for all $i\geq 1$. Given a
tower of function fields we write $N^{(i)}=N(F^{(i)})$,
$g^{(i)}=g(F^{(i)})$ and we say that the tower is asymptotically good
if $g^{(i)} \rightarrow \infty$ for $i \rightarrow \infty$ and 
$\liminf_{i\rightarrow \infty}(N^{(i)}/g^{(i)})=\kappa$ holds for some
$\kappa > 0$. We mention that the interest in asymptotically good towers of 
function fields partly comes from the fact that they give rise to
arbitrary long codes with good parameters (see~\cite[Sec.\
VII.2]{stichtenoth} for the details). We now present a corollary to Theorem~\ref{themain} concerning asymptotically good towers.
\begin{Corollary}\label{corfour}
Assume a tower of  function fields is given with $g^{(i)}
\rightarrow \infty$ for $i \rightarrow \infty$ and  $\liminf_{i
  \rightarrow \infty}(\frac{N^{(i)}}{g^{(i)}}) = \kappa > 0$ (that is,
the tower is asymptotically good). Let $(\PP^{(1)},\PP^{(2)}, \ldots )$ be
any sequence such that $\PP^{(i)}$ is a rational place of $F^{(i)}$ for
$i=1,2, \ldots $. Let $\lambda_1^{(i)}$ be the multiplicity of the
Weierstrass semigroup related to $\PP^{(i)}$ and let $m_i$ be the number
of generators in some description of $\Lambda^{(i)}$. We have
\begin{eqnarray}
&\liminf_{i \rightarrow \infty}(\lambda_1^{(i)}/g^{(i)}) \geq \kappa /q
\label{eqsnabel3}\\
&m_i \rightarrow \infty {\mbox{ \ for \  }} i \rightarrow \infty \label{eqmiuendelig}
\end{eqnarray} 
\end{Corollary}
\noindent {\it{Proof:}} \ \ From Lewittes' bound~(\ref{stjerneskud3}) we know that $N^{(i)} \leq q
\lambda_1^{(i)}+1$. Applying the assumptions from the corollary we get
$$
\liminf_{i \rightarrow \infty}(\frac{q\lambda_1^{(i)}+1}{g^{(i)}})
\geq \kappa \Rightarrow \liminf_{i \rightarrow \infty}(\frac{\lambda_1^{(i)}}{g^{(i)}})
\geq \frac{\kappa}{q}. 
$$
To see the last part of the corollary observe that by
Proposition~\ref{proboundbound} $N^{(i)} \leq q^{m_i}+1$ holds.
\ \firkant
We next apply Corollary~\ref{corfour} to the case of the asymptotically good tower of function fields over $\lmto$ which was introduced by Garcia
and Stichtenoth in~\cite{GS2}. This tower is defined by
$F^{(1)}=\lmto(x_1)$ and $F^{(i+1)}=F^{(i)}(x_{i+1})$ with 
$$x_{i+1}^q+x_{i+1}=\frac{x_i^q}{x_i^{q-1}+1}.$$
The tower is actually as good as a tower over $\lmto$ can possibly
be. We have 
\begin{equation}
g^{(i)}=\left\{ \begin{array}{ll}
(q^{i/2}-1)^2 & {\mbox{ \ if\ }} i \equiv 0 {\mbox{\ mod \ }} 2 \\
(q^{(i+1)/2}-1)(q^{(i-1)/2}-1) & {\mbox{ \ if\ }} i \equiv 1 {\mbox{\
    mod \ }} 2 
\end{array} \right.
\label{genus}
\end{equation}
which of course implies $g^{(i)} \rightarrow \infty$ for $i \rightarrow \infty$ and we have 
\begin{equation}
\lim_{i \rightarrow \infty}(N^{(i)}/g^{(i)})=q-1. \label{factor}
\end{equation}
The element $x_1 \in F^{(1)} \subseteq F^{(i)}$ has in $F^{(i)}$ a
unique pole which we denote by $\PP_{\infty}^{(i)}$. This place is known
to be rational. Our interest in $\PP_{\infty}^{(i)}$ comes from the convenient
fact that the Weierstrass semigroup $\Lambda^{(i)}$ of
$\PP_{\infty}^{(i)}$ was established in~\cite{PST}. This allows us to
apply Corollary~\ref{corfour} and Remark~\ref{remone}. For $i \geq 1$
define
\begin{eqnarray}
c^{(i)}&=&\left\{ \begin{array}{ll}
q^i-q^{i/2}& {\mbox{\ if \ }} i \equiv 0 {\mbox{ \ mod \ }} 2 \\
q^i-q^{(i+1)/2}& {\mbox{\ if \ }} i \equiv 1 {\mbox{ \ mod \ }} 2
\end{array} \right. \nonumber
\end{eqnarray}
then $\Lambda^{(1)}={\mathbb{N}}_0$, and for $i \geq 1$ we have
$$\Lambda^{(i+1)}= q\Lambda^{(i)} \cup \{ x \in {\mathbb{N}}_0 \mid x
\geq c^{(i+1)} \}.$$
Assume for a
moment that $q >2$. The three smallest non gaps of $\Lambda^{(i)}$ are
$0, q^{i-1}, 2q^{i-1}$ (this can be verified by an induction proof).
Hence, if we write
 $\Lambda^{(i)}= \langle \lambda_1^{(i)}, \ldots ,
\lambda_{m(i)}^{(i)} \rangle$ with $0< \lambda_1^{(i)} < \cdots <
\lambda_{m(i)}^{(i)}$ then we get $\lambda_1^{(i)}=q^{i-1}$ and
$\lambda_2^{(i)} \geq 2g^{i-1}$. 
From~(\ref{genus}) we see that
$g^{(i)} < q^i$ holds. Therefore $$q^2 \lambda_1^{(i)}+2g^{(i)} \leq
q^2\lambda_2^{(i)},$$
and by Remark~\ref{remone} the bounds~(\ref{stjerneskud2}) and
(\ref{stjerneskud3}) from Theorem~\ref{themain} will produce the same
results when applied to $\PP_{\infty}^{(i)}$ in the cases $q>2$. 
Hence, there is no point in trying to improve upon
Corollary~\ref{corfour} for the case of Garcia and Stichtenoth's
second tower by taking into account the new
bound~(\ref{stjerneskud2}).\\
Whether or not $q >2$ or $q=2$ holds we have
$\lambda_1^{(i)}=q^{i-1}$, and therefore
$$\lim_{i \rightarrow \infty}(\lambda_1^{(i)}/g^{(i)}) = 1/q.$$
For comparison Corollary~\ref{corfour} reads
$$\liminf_{i \rightarrow \infty}(\lambda_1^{(i)}/g^{(i)}) \geq
(q-1)/q^2=\frac{1-1/q}{q}.$$ 
For large values of $q$ Corollary~\ref{corfour} therefore gives a
reasonable picture of the situation in the case of the tower considered.\\
We conclude this section by showing that one can not hope to produce
asymptotically good towers having telescopic Weierstrass
semigroups.  This fact may be known to some of the
researchers of asymptotically good towers; but we have not been able
to find any reference on it.
\begin{Definition}\label{deftele} Let $\Lambda=\langle a_1, \ldots , a_k \rangle \in
{\mathbb{N}}_0$ be a semigroup for which we have 
 $\gcd(a_1, \ldots , a_k)=1$. For $1 \leq j \leq k$ define $d_j:=\gcd (a_1, \ldots , a_j)$
and $\Lambda_j:=\langle a_1/d_j, \ldots , a_j/d_j \rangle$. If
$a_j/d_j \in \Lambda_{j-1}$ for $2 \leq j \leq k$ then $\Lambda$ is
said to be telescopic.
\end{Definition}
We will need the following result corresponding to~\cite[Lem.\
5.34]{handbook}.
\begin{Lemma}\label{lemtele}
If $\Lambda=\langle a_1, \ldots , a_k \rangle$ describes a telescopic
semigroup as in Definition~\ref{deftele} then for any $\lambda \in \Lambda$
there exist (uniquely) determined non negative integers $x_1, \ldots ,
x_k$ such that $0 \leq x_j \leq d_{j-1}/d_j$ for $2 \leq j \leq k$ and
$\lambda= \sum_{j=1}^kx_ja_j$. 
\end{Lemma}
\begin{Proposition}\label{protele}
Let $(F^{(1)}/{\mathbb{F}}_q,F^{(2)}/{\mathbb{F}}_q,\ldots )$ be a
tower of function fields such that for all $i$ (or
alternatively at least infinitely many $i$) the following
holds: $F^{(i)}$ possesses a rational place $\PP^{(i)}$ having a
telescopic Weierstrass semigroup $\Lambda^{(i)}$. Then the tower is
asymptotically bad. 
\end{Proposition}
\noindent {\it {Proof:}} \ \ Only in the case that $N^{(i)}
\rightarrow \infty$ for $i \rightarrow \infty$ we can hope to get an
asymptotically good tower. We therefore assume that $N^{(i)}$
satisfies this condition. As $\Lambda^{(i)}$ is telescopic we know
that there exists a description $\Lambda^{(i)}= \langle
\lambda_1^{(i)}, \ldots  , \lambda_{m_i}^{(i)} \rangle$ satisfying the
conditions in Definition~\ref{deftele}. We will assume that we have
chosen a description of this kind with $m_i$ smallest
possible. Following Definition~\ref{deftele} we have $d_j^{(i)}:=\gcd
(\lambda_1^{(i)}, \ldots , \lambda_{j}^{(i)} )$ for $1 \leq j \leq
m_j$. Clearly, $d_j^{(i)} \mid d_{j-1}^{(i)}$ for $j \geq 2$ and by
minimality of $m_i$ Lemma~\ref{lemtele} implies $d_{j-1}^{(i)} \geq 2
d_j^{(i)}$. The genus $g^{(i)}$ is given by the following closed form
expression (see~\cite[Pro.\ 5.35]{handbook})
$$g^{(i)}=\left( 1 +
  \sum_{j=2}^{m_i}(\frac{d_{j-1}}{d_j}-1)\lambda_j^{(i)} \right) /2,$$
and therefore $g^{(i)} \geq \frac{m_i-1}{2}\lambda_1^{(i)}$ holds. On
the other hand Lewittes' bound~(\ref{stjerneskud3}) states $N^{(i)}
\leq q \lambda_1^{(i)}+1$. Hence, $\lim_{i \rightarrow
  \infty}\frac{N^{(i)}}{g^{(i)}}=0$ 
follows immediately from (\ref{eqmiuendelig}).\\
\ \ \ \ \  \firkant 
\section{Concluding remarks and acknowledgements}\label{seccon}\label{secon}
We would like to mention that it is possible to give an alternative proof of Theorem~\ref{themain}
by using results on order domain theory from~\cite{GP} in combination
with the footprint bound from Gr\"{o}bner basis theory. It would be nice if the methods from the present paper can to some extend be
used to deal with algebraic function fields of more variables. We
leave it for further research to decide if this is possible. The
gonality of a curve is closely related to the multiplicities
$\lambda_1$ of the Weierstrass semigroups studied in the present paper. More generally the notion of
the gonality sequences that one finds in the papers~\cite{munuera1},
\cite{munuera2} and~\cite{yang} on generalized Hamming Weights is
closely related to the generators $\lambda_1, \ldots , \lambda_m$
of the Weierstrass semigroups studied in the present paper. We leave it for further research to
decide if the new bound~(\ref{stjerneskud2}) has some implications for
the theory of gonality sequences.\\ 
The
authors would like to thank Peter Beelen, Tom H{\o}holdt, Massimiliano
Sala and Ruud Pellikaan for pleasent discussions.


\begin{thebibliography}{99}
\bibitem{GS2} A.\ Garcia and H.\ Stichtenoth, ``On the asymptotic
  behaviour of some towers of function fields over finite fields,''
  {\em{J.\ Number Theory}}, {\textbf{61}} (1996), pp.\ 248--273.
\bibitem{geervlugt} G.\ van der Geer and M.\ van der Vlugt, ``Tables
  of Curves with Many Points,''\\
\url{http://www.science.uva.nl/~geer/tables-mathcomp18.pdf} \\
(August 20, 2007).
\bibitem{normtrace} O.\ Geil, ``On codes from norm-trace curves,''
  {\em{Finite Fields Appl.,}} {\textbf{9}} (2003), pp.\ 351--371.
\bibitem{GP} O.\ Geil and R.\ Pellikaan, ``On the structure of order
  domains,'' {\em{Finite Fields Appl.,}} {\bf{8}} (2002), pp.\ 369--396.
\bibitem{handbook} T.\ H{\o}holdt, J.\ van Lint and R.\  Pellikaan, Algebraic Geometry Codes, Chapter 10 in 
    ``Handbook of Coding Theory,'' (V.S.\ Pless and W.C.\ Huffman,
  Eds.), vol.\ 1, Elsevier, Amsterdam, 1998, pp.\ 871--961.
\bibitem{lewittes} J.\ Lewittes, ``Places of degree one in function
  fields over finite fields,'' {\em{J.\ Pure Appl.\ Algebra}},
      {\textbf{69}} (1990), pp.\ 177--183.
\bibitem{munuera1} C.\ Munuera, ``On the generalized Hamming weights of
        geometric Goppa codes,'' {\em{IEEE Trans.\ Inform.\ Theory,}}
    {\textbf{45}} (1994), pp.\ 2092--2099.
\bibitem{munuera2} C.\ Munuera and D.\ Ramirez, ``The second
      and third generalized Hamming weights of Hermitian codes,''
  {\em{IEEE Trans.\ Inform.\ Theory,}} {\textbf{45}} (1999), pp.\ 709--712.
\bibitem{oesterleserre} H.\ Niederreiter and C.\ Xing,
  {\textit{Rational Points on Curves over Finite Fields - Theory and
      Applications}}, London Math.\  Soc.\  Lecture Note Series,
  {\textbf{285}}, Cambridge University Press, 2001.
\bibitem{rivaldo} Nivaldo Mediros, homepage:\\
\url{http://w3.impa.br/~nivaldo/algebra/semigroups/} \\
(October 1, 2007).
\bibitem{PST} R.\ Pellikaan, H.\ Stichtenoth and F.\ Torres,
  ``Weierstrass semigroups in an asymptotically good tower of function
  fields,'' {\em{Finite Fields Appl.}}, {\textbf{4}}
  (1998), pp.\ 381--392.
\bibitem{stichtenoth} H.\ Stichtenoth, {\textit{Algebraic Function
      Fields and Codes}}, Springer Verlag, 1993.
\bibitem{storvoloch} K.\ O.\ St\"{o}hr and J.\ F.\ Voloch, ``Weierstrass
  points and curves over finite fields,'' {\em {Proc. London Math.\
      Soc.}} (3), {\textbf{52}} (1986), pp.\ 1--19.
\bibitem{yang} K.\ Yang, P.\ V.\ Kumar and H.\ Stichtenoth,
  ``On the weight hierarchy of geometric Goppd codes,''
  {\em{IEEE Trans.\ Inform.\ Theory,}} {\textbf{40}} (1994), pp.\ 913--920.
\end{thebibliography}
\end{document}